\documentclass[reqno]{amsart}

\pdfoutput=1

\usepackage[alphabetic]{amsrefs}
\usepackage{amsmath}
\usepackage{amsfonts}
\usepackage{amssymb}
\usepackage{enumerate}
\usepackage{amstext}
\usepackage{amsbsy}
\usepackage{amsopn}
\usepackage{bbm,amsthm}
\usepackage{amscd}
\usepackage[pdftex]{color}
\usepackage[pdftex]{graphicx}
\usepackage[all]{xy}
\usepackage{mathtools}

\usepackage{tikz}
\usetikzlibrary{intersections,positioning,patterns,arrows,decorations.markings}

\usepackage{todonotes}

\usepackage{amsxtra}
\usepackage{upref}
\usepackage{epstopdf}
\usepackage{graphicx,color}
\usepackage{hyperref}
\usepackage{longtable}
\usepackage[francais, english]{babel}

\makeatletter
\newsavebox{\@brx}
\newcommand{\llangle}[1][]{\savebox{\@brx}{\(\m@th{#1\langle}\)}%
  \mathopen{\copy\@brx\kern-0.5\wd\@brx\usebox{\@brx}}}
\newcommand{\rrangle}[1][]{\savebox{\@brx}{\(\m@th{#1\rangle}\)}%
  \mathclose{\copy\@brx\kern-0.5\wd\@brx\usebox{\@brx}}}
\makeatother

\DeclareFontFamily{OML}{rsfs}{\skewchar\font'177}
\DeclareFontShape{OML}{rsfs}{m}{n}{ <5> <6> rsfs5 <7> <8> <9> rsfs7
  <10> <10.95> <12> <14.4> <17.28> <20.74> <24.88> rsfs10 }{}
\DeclareMathAlphabet{\mathfs}{OML}{rsfs}{m}{n}

\newtheorem{theorem}{Theorem}

\theoremstyle{definition}
\newtheorem{definition}[theorem]{Definition}

\theoremstyle{remark}

\numberwithin{equation}{section}
\numberwithin{theorem}{section}

\newcommand{\intav}[1]{\mathchoice {\mathop{\vrule width 6pt height 3 pt depth  -2.5pt
\kern -8pt \intop}\nolimits_{\kern -6pt#1}} {\mathop{\vrule width
5pt height 3  pt depth -2.6pt \kern -6pt \intop}\nolimits_{#1}}
{\mathop{\vrule width 5pt height 3 pt depth -2.6pt \kern -6pt
\intop}\nolimits_{#1}} {\mathop{\vrule width 5pt height 3 pt depth
-2.6pt \kern -6pt \intop}\nolimits_{#1}}}

\newcommand{\intavl}[1]{\mathchoice {\mathop{\vrule width 6pt height 3 pt depth  -2.5pt
\kern -8pt \intop}\limits_{\kern -6pt#1}} {\mathop{\vrule width 5pt
height 3  pt depth -2.6pt \kern -6pt \intop}\nolimits_{#1}}
{\mathop{\vrule width 5pt height 3 pt depth -2.6pt \kern -6pt
\intop}\nolimits_{#1}} {\mathop{\vrule width 5pt height 3 pt depth
-2.6pt \kern -6pt \intop}\nolimits_{#1}}}

\newcommand{\vertiii}[1]{{\left\vert\kern-0.2ex\left\vert\kern-0.2ex\left\vert #1 
    \right\vert\kern-0.2ex\right\vert\kern-0.2ex\right\vert}}

\newcommand{\ve}{\varepsilon}

\newcommand{\R}{\mathbb{R}}

\newcommand{\hd}{{\rm HD}}
\newcommand{\te}{\theta}
\newcommand{\proj}{\text{proj}}

\begin{document}

\title[$d$--sets with irregular projections of Hausdorff measures
]{Examples of $d$--sets with irregular\\ projection of Hausdorff measures}

\author{Yuri Lima and Carlos Gustavo Moreira}

\address{Departamento de Matem\'atica, Universidade Federal do Cear\'a (UFC), Campus do Pici,
Bloco 914, CEP 60455-760. Fortaleza -- CE, Brasil}
\email{yurilima@gmail.com}
\address{IMPA, Estrada Dona Castorina 110, CEP 22460-320, Rio de Janeiro, Brazil}
\email{gugu@impa.br}

\date{\today}
\thanks{Lima is supported by CNPq and Instituto Serrapilheira, grant ``Jangada Din\^amica:
Impulsionando Sistemas Din\^amicos na Regi\~ao Nordeste''. Moreira is supported by CNPq, FAPERJ,
and INCTMAT project of J. Palis}

\begin{abstract}
Given positive integers $\ell<n$ and a real $d\in (\ell,n)$, we construct
sets $K\subset \R^n$ with positive and finite Hausdorff $d$--measure such that the Radon-Nikodym derivative
associated to all projections on $\ell$--dimensional planes is not an $L^p$ function, for
all $p>1$.
\end{abstract}

\maketitle

\section{Introduction}\label{Section-introduction}

If $U$ is a subset of $\R^n$, the diameter of $U$ is $|U|=\sup\{|x-y|;x,y\in U\}$
and, if $\mathcal U$ is a family of subsets of $\R^n$, the {\it diameter} of $\mathcal U$ is defined as
$$\left\|\mathcal U\right\|=\sup_{U\in\,\mathcal U}|U|.$$
Given $d>0$, the {\it Hausdorff $d$--measure} of a set $K\subset\R^n$ is
$$m_d(K)=\lim_{\ve\rightarrow 0}\left(\inf_{\mathcal U\text{ covers }K\atop{\left\|\mathcal U\right\|<\,\ve}}
\sum_{U\in\,\mathcal U}|U|^d\right).$$
In particular, when $n=1$, $m=m_1$ is the Lebesgue measure of Lebesgue measurable sets on $\R$.
It is not difficult to show that there exists a unique $d_0\ge0$ for which $m_d(K)=\infty$ if $d<d_0$ and
$m_d(K)=0$ if $d>d_0$. We define the Hausdorff dimension of $K$ as $\hd(K)=d_0$. For each $\te\in\R$,
let $v_\te=(\cos \te,\sin\te)$, $L_\te$ be the line in $\R^2$ through the origin containing $v_\te$ and
$\proj_\te:\R^2\rightarrow L_\te$ be the orthogonal projection. Noting that $L_{\te}=L_{\te+\pi}$, from now on
we restrict $\te$ to the interval $[-\pi/2,\pi/2]$.

In $1954$, J. M. Marstrand \cite{Marstrand} proved the following result on the fractal dimension of plane sets.

\begin{theorem}[Marstrand]\label{thm 1}
If $K\subset\R^2$ is a Borel set with $\hd(K)>1$, then
$m(\proj_\te(K))>0$ for $m$--almost every $\te\in\R$.
\end{theorem}

In the article \cite{Lima-Moreira}, we gave a new proof of Theorem \ref{thm 1}, using elementary combinatorial tools.
Since this theorem is about lower bounds for Hausdorff dimension and Lebesgue measure,
it is sufficient to prove it for some suitable subset of $K$. We used this reduction, restricting to a subset
$K'\subset K$ that is more regular, which we called $d$--regular, when:
\begin{enumerate}[i.]
\item[(i)] $m_d(K')>0$;
\item[(ii)] there exists $b>0$ such that $m_d(K'\cap B(x,r))\leq br^d$ for all $x\in K$ and $r>0$, where $B(x,r)\subset\R^2$
is the open ball with center $x$ and radius $r$. 
\end{enumerate}
Under these assumptions, we can apply a double counting argument and then prove Theorem \ref{thm 1}.

The restriction of $m_d$ to $K$ induces by $\proj_\te$ a measure 
$\mu_\te:=(\proj_\te)_*(m_d\restriction_K)$ on $L_\te$. In Theorem 1.2 of \cite{Lima-Moreira}, we proved
that $\mu_\te$ is absolutely continuous with respect to $m$ for $m$--almost every $\theta$.
We also {\em erroneously} claimed
that the Radon-Nikodym derivative $f_\theta=\tfrac{d\mu_\te}{dm}$ is an $L^2$ function
for $m$--almost every $\te$. As a matter of fact, we proved this property for the $d$--regular subset $K'\subset K$, 
but this {\em does not} imply the same result for the larger set $K$.

The goal of this paper is to show that, indeed, there are examples of sets $K$ for which $f_\theta$
is not an $L^2$ function for any $\te$. We actually construct examples in arbitrary dimension
such that the respective Radon-Nikodym derivative, when it exists, is not an $L^p$ function for any $p>1$.
Before stating the theorem, we need some notation.

\begin{definition}[$d$--set]
A Borel set $K\subset \R^n$ is called a {\em $d$--set} if $\hd(K)=d$ and $0<m_d(K)<\infty$.
\end{definition}

Given positive integers $\ell<n$, let $G(n,\ell)$ be the Grassmannian of $\ell$--dimensional
planes in $\R^n$ containing the origin. Given $\pi\in G(n,\ell)$, we write $\Pi:\R^n\to\pi$ for the orthogonal projection
onto $\pi$ and $m_\ell$ for the associated $\ell$--dimensional Lebesgue measure on $\pi$. 
Let $\ell<d<n$.

\begin{theorem}[Main theorem]\label{thm-main}
Let $(c_j),(r_j)$ be sequences of positive real numbers such that $\sum c_jr_j^d<+\infty$
and $c_j r_j^{d-\ve}\to +\infty$, for all $\ve>0$.
Let $K=\bigcup K_j$ be a disjoint union of sets $K_j$, where $K_j\subset \R^n$ is a $d$--set with 
${\rm diam}(K_j)\leq r_j$ and
$m_d(K_j)=c_jr_j^d$, and let $\mu_\pi=\Pi_*(m_d\restriction_K)$.
Then $K$ is a $d$--set with $0<m_d(K)<+\infty$ such that for every $\pi\in G(n,\ell)$
the Radon-Nikodym derivative $d\mu_\pi/dm_\ell$, when it exists, is not in $L^p$ 
for any $p>1$. 
\end{theorem}

The proof of Theorem \ref{thm-main} is contained in the next section.
By the above assumptions, necessarily $r_j\to 0$ and $c_j\to+\infty$, so $K$ is a set that does not
satisfy assumption (ii). We can take e.g. $r_j=2^{-j}$ and $c_j=2^{jd}/j^2$, 
in which case $c_jr_j^d=j^{-2}$ and $c_j r_j^{d-\ve}=2^{\ve j}/j^2$.

For the sake of completeness, let us show that, given a positive integer $n$ and any positive
numbers $r, c, d$ with $0<d<n$, there is a $d$--set $Z\subset\mathbb R^n$ with ${\rm diam}(Z)\le r$ and $m_d(Z)=c$.
Start with a compact $d$--subset $X$ of $(0,1)^n$ with $m_d(X)=a>0$,
say a Cartesian product of $n$ homogeneous regular Cantor sets contained in $(0,1)$
with Hausdorff dimensions $d/n$.
Take a cube of side $r/\sqrt{n}$ (and so with diameter $r$), divide it in $M^n$ cubes of sides $r/M\sqrt{n}$
(where $M$ is a large integer to be chosen later), and put a homothetic copy (with ratio $r/M\sqrt{n}$) 
of $X$ in each of these small cubes. The resulting set $Y$ has 
$m_d(Y)=M^na(r/M\sqrt{n})^d=M^{n-d} a r^d n^{-d/2}$,
which can be made larger than $c$ by taking $M$ large enough.
Then a contractive homothety sends $Y$ to a $d$--set $Z$ with diameter smaller than $r$ satisfying $m_d(Z)=c$.

We also take the chance to state Theorem 1.2 of \cite{Lima-Moreira} correctly.

\begin{theorem}[Theorem 1.2 of \cite{Lima-Moreira} corrected]
Let $K\subset \R^2$ be a Borel set with ${\rm HD}(K)>1$, and let $\mu_\theta$ as above.
Then $\mu_\theta$ is absolutely continuous with respect to $m$, for $m$--almost every $\theta\in\R$.
If furthermore
there exists $b>0$ such that $m_d(K\cap B(x,r))\leq br^d$ for all $x\in K$ and $r>0$,
then the Radon-Nikodym derivative $d\mu_\te/dm$ is an $L^2$ function, for
$m$--almost every $\te\in\R$.
\end{theorem}

\section{Proof of Theorem \ref{thm-main}}
It is clear that $m_d(K)\in (0,+\infty)$, since
$$
0<m_d(K_1)\leq m_d(K)\leq \sum m_d(K_j)=\sum c_jr_j^d<+\infty.
$$
For the second part, fix $\pi\in G(n,\ell)$, $p>1$, and write $f=d\mu_\pi/dm_\ell$ 
(assuming that $\mu_\pi$ is absolutely continuous with respect to Lebesgue).
Let $\mu_{\pi,j}=\Pi_*(m_d\restriction_{K_j})$,
and let $f_j=d\mu_{\pi,j}/dm_\ell$ (notice that $\mu_{\pi,j}\ll \mu_\pi$). Each $f_j\geq 0$ and $f=\sum f_j$. 
We will obtain an estimate
for the $L^p$ norm of $f_j$. First, observe that
$$
\|f_j\|_1=\int_{\Pi(K_j)} \frac{d\mu_{\pi,j}}{dm_\ell}dm_\ell=\mu_\Pi[\Pi(K_j)]=m_d(K_j).
$$
Now, fix a ball $B_j\supset K_j$ of radius $r_j$, and let $g=1_{\Pi(B_j)}$.
Since $\Pi(B_j)$ is a ball of radius $r_j$, for any $q\geq 1$ we have
$$
\|g\|_q^q=\|g\|_1=m_\ell[\pi(B_j)] \leq Cr_j^\ell,
$$
where $C>0$ is a constant that only depends on $\ell$. 
Taking $q>1$ such that $\tfrac{1}{p}+\tfrac{1}{q}=1$, the H\"older inequality 
implies that $\|f_j\|_p \|g\|_q \geq \|f_j\|_1$ and so
\begin{eqnarray*}
\|f_j\|_p\geq \frac{\|f_j\|_1}{\|g\|_q}\geq \frac{m_d(K_j)}{C^{\frac{1}{q}}r_j^{\frac{\ell}{q}}}
=C^{-\frac{1}{q}}c_jr_j^{d-\frac{\ell}{q}}.
\end{eqnarray*}
Since $f_j\geq 0$, we conclude by the standing assumption with $\ve=\tfrac{\ell}{q}$ that
$$
\|f\|_p\geq \limsup_{j\to+\infty}\|f_j\|_p\geq  \lim_{j\to+\infty}C^{-\frac{1}{q}}c_jr_j^{d-\frac{\ell}{q}}=
C^{-\frac{1}{q}}\lim_{j\to+\infty}c_jr_j^{d-\frac{\ell}{q}}=+\infty.
$$
This concludes the proof of the theorem.

\bibliographystyle{alpha}
\bibliography{bibliography}{}

\end{document}